\documentclass{amsart}

\usepackage{latexsym}
\usepackage{amsmath}
\usepackage{amssymb}

\usepackage{multicol}
\usepackage{makecell}

\usepackage[cp1251]{inputenc}
\usepackage[english]{babel}

\theoremstyle{plain}

\newtheorem{theorem}{Theorem}[section]
\newtheorem{lem}[theorem]{Lemma}
\newtheorem{pr}[theorem]{Proposition}

\newtheorem*{theorem*}{Theorem}
\newtheorem*{lem*}{Lemma}

\newtheorem{rem}[theorem]{Remark}

\numberwithin{equation}{section}

\newcommand{\er}{\varepsilon}
\newcommand{\si}{\sigma}
\newcommand{\sid}{\sigma_d}

\newcommand{\za}{\zeta}

\newcommand{\ph}{\varphi}
\newcommand{\hol}{\mathcal{H}ol}
\newcommand{\Nbb}{\mathbb N}
\newcommand{\Rbb}{\mathbb R}
\newcommand{\Dbb}{\mathbb D}

\newcommand{\T}{\mathbb{T}}

\newcommand{\bd}{B_d}
\newcommand{\spd}{\partial B_d}
\newcommand{\ccd}{\mathbb{C}^d}
\newcommand{\bbd}{\overline{B}_d}
\newcommand{\MM}{M_+}

\newcommand{\hp}{{H^p}}

\newcommand{\rad}{\mathcal{R}}

\newcommand{\KK}{\mathcal{K}}
\newcommand{\bdd}{\overline{B}_d}

\begin{document}

\date{}


\author{Evgueni Doubtsov}
\address{St.~Petersburg Department
of Steklov Mathematical Institute, Fontanka 27, 
St.~Petersburg 191023, Russia}
\email{dubtsov@pdmi.ras.ru}


\title[Reverse Carleson measures]{Reverse Carleson measures in several variables}

\begin{abstract}
Let $q > 0$ and
let $X$ be a quasi-Banach space of analytic functions in the unit ball $B_d$ of $\mathbb{C}^d$.
A finite positive Borel measure $\mu$ on the closed unit ball $\overline{B}_d$
is called a \emph{$q$-reverse Carleson measure} for $X$ if there exists a constant $C > 0$ such that
\[
\|f\|_{X}\leq C \|f\|_{L^q(\overline{B}_d, d\mu)}
\]
for all $f\in X\cap C({\overline{B}_d})$.
For $0<q< \infty$, we characterize the $q$-Carleson measures for the Hardy spaces $H^p(B_d)$, $0<p<\infty$.
\end{abstract}


\maketitle

\section{Introduction}\label{s_int}
	
Let $\bd$ denote the open unit ball of $\ccd$, $d\ge 1$,
and let $\si=\sid$ denote the normalized Lebesgue measure on the unit sphere $\spd$. We denote by $\hol(\bd)$ the space of holomorphic functions in the ball $\bd$.
Also, put $\Dbb = B_1$ and $\T = \partial \Dbb$.

\subsection{Hardy spaces and direct Carleson measures}
For $0<p<\infty$ the Hardy space $H^p=H^p(\bd)$
consists of functions $f\in \hol(\bd)$ satisfying
\[
\|f\|_{H^p}^p = \sup_{0<r<1} \int_{\spd} |f(r\zeta)|^p\, d\si(\zeta) < \infty.
\]
As usual we identify the Hardy space $H^p(\bd)$, $p>0$, with the space
$H^p(\spd)$ which consists of the corresponding boundary values.
We refer the reader to the monographs \cite{Ru80} and \cite{Z05} for more information on the Hardy spaces $\hp(\bd)$.

For an arbitrary dimension $d\ge 1$, the following sets play a role similar to that of  
the subarcs of the unit circle $\T$:
\begin{equation}\label{e_qball}
Q = Q(\za, \delta) =\{\xi\in \spd: |1-\langle \za, \xi\rangle | \le \delta\},
\quad \za\in\spd,\ \delta>0.
\end{equation}
Remark that $Q$ is a closed ball with respect to the non-isotropic metric
\[
\rho(\za, \xi) = |1-\langle \za, \xi\rangle |^{\frac{1}{2}}, \quad \za, \xi\in\spd.
\]
Consider $\MM(\bd)$, the space of finite, positive Borel measures on the ball $\bd$.
The theorem about Carleson measures in the unit ball (see for instance \cite[Chapter~5]{Z05})
yields that the bound
\begin{equation}\label{e_Carl_n}
\| f \|_{L^p(B_d, \mu)} \le C \|f\|_\hp,\quad f \in\hp,
\end{equation}
is equivalent to the following condition on the measure $\mu$:
\begin{equation}\label{e_Carl_geom_n}
\mu(S_Q) \le C\si(Q)\quad\textrm{for all non-isotropic balls}\ Q\subset\spd,
\end{equation}
where
\[
S_Q = \left\{z \in\bd:\ 1-\delta \le |z| < 1,\ \frac{z}{|z|} \in  Q \right\}
\]
is a standard Carleson window for $Q= Q(\cdot, \delta)$.
Remark that  for $d=1$ the above result was obtained by L.~Carleson \cite{Ca58, Ca62}.

By now, the Carleson measures are characterized for various spaces 
on the unit disc and other domains.

\subsection{Hardy spaces and reverse Carleson measures}
In the present paper, we are interested in the following reverse Carleson inequality:
\begin{equation}\label{e_reverse}
\|f \|_\hp \le C \|f\|_{L^q(\bbd,\mu)}, \quad f\in  C(\bbd) \cap H^p(\bd),\ 0 < p, q < \infty,
\end{equation}
for $\mu\in\MM(\bbd)$.
Note that $\mu$ is defined on the closed unit ball $\bbd$.
The measures with property \eqref{e_reverse} were characterized in \cite{HMNO14} for $d=1$ and $p=q>1$.
The case $d\ge 1$ and $p=q>0$ was recently considered in \cite{Dou24} and \cite{DTV26}.
In particular, \eqref{e_reverse} with $p=q>0$ is equivalent to the following property:
\[
\mu(Q) \ge C\si(Q)
\]
for an absolute constant $C$ and for all non-isotropic balls $Q\subset\spd$.

In Section~\ref{s_pq}, we describe the $q$-reverse Carleson measures for the space $H^p(\bd)$. 
We obtain a complete answer to this problem for all values $0 < p, q < \infty$. The argument is based on
certain ideas from harmonic analysis.
For the results concerning the corresponding problem for direct Carleson measures, see \cite{Dur69, Lu91}.

 	\subsection*{Notation}
	In what follows, $C>0$ denotes an absolute constant, the value of which may vary from line to line.

	\section{$q$-reverse Carleson measures for Hardy spaces}\label{s_pq}
	
The main goal of the present section is to describe the $q$-reverse Carleson measures for the spaces $H^p(\bd)$.
As mentioned above, it suffices to consider the case $p\neq q$.
	
	\subsection{The balayage phenomenon}

We need the following version of Ryll--Wojtasz\-czyk polynomials. 

\begin{pr}[see e.g. \cite{A86}]\label{p_RW}
Let $d\ge 2$.
There exist $J \in \Nbb$ and $\delta>0$ such that for each $n\in\Nbb$
there exist holomorphic homogeneous polynomials $W_{n j}$, $1\le j \le J$,
of degree $n$ on $\ccd$ and with following properties:
\begin{align}
\sum_{j=1}^J |W_{n j}(\za)| &\le 1, \quad \za\in \spd, \label{e_RWup}\\
\max_{1\le j \le J} |W_{n j}(\za)| &\ge \delta>0,   \quad \za\in \spd.  \label{e_RWlow}
\end{align}
\end{pr}

	We start with the following observation: 
	only the restriction $\mu_{|\spd}$ affects the inequality \eqref{e_reverse}.
	
	\begin{lem}
		\label{general_balayage}
		Suppose that for a measure $\mu\in M_+(\bdd)$ there exists a constant $C > 0$ such that 
		\begin{equation}\label{q-RCM-Hardy}
			\|f\|_{H^p} \le C \|f\|_{L^q(\mu)}
		\end{equation}
		for every $f\in H^p(\bd)\cap C(\bdd)$.
		Put $\nu=\mu|_{\spd}$. Then the inequality
		\begin{equation}
			\label{q-RCM-Hardy-balayaged}
			\|f\|_{H^p} \le C\|f\|_{L^q(\nu)}
		\end{equation}
		holds for every $f\in H^p(\bd)\cap C(\bdd)$.
	\end{lem}
	
	\begin{proof}
		Let $f\in H^p(\bd)\cap C(\bdd)$. 
Consider the following functions: 
\[
f_{n j}(z)= W_{n j}(z) f(z), \quad j=1, \dots, J,
\] 
where $W_{n j}$ are provided by Proposition~\ref{p_RW}.
Applying the inequality \eqref{q-RCM-Hardy} to the functions $f_{n j}$, we obtain
		\begin{equation}\label{RCM_fN}
			\|f_{n j}\|_{H^p} \le C \|f_{n j}\|_{L^q(\mu)},
		\end{equation}
hence,
\begin{equation}\label{RCM_sum}
			\sum_{j=1}^J \|f_{n j}\|_{H^p} \le C \sum_{j=1}^J \|f_{n j}\|_{L^q(\mu)}.
		\end{equation}

Now, on the one hand, observe that
\begin{equation}\label{e_left}
\sum_{j=1}^J \|f_{n j}\|_{H^p} \ge C \left( \int_{\spd} \sum_{j=1}^J |W_{n j}(z) f(z)|^p\, d\si \right)^{\frac{1}{p}}
 \ge \delta C \left( \int_{\spd} |f(z)|^p\, d\si \right)^{\frac{1}{p}}
\end{equation}
by \eqref{e_RWlow}.

On the other hand, by \eqref{e_RWup},
\begin{equation}\label{e_right}
 \sum_{j=1}^J \|f_{n j}\|_{L^q(\nu)} 
\le C \|f\|_{L^q(\nu)}.
\end{equation}
Next, by \eqref{e_RWup}, we have for every $z\in \bd$ and $j=1, \dots, J$,
\[
|f_{n j}(z)|=|W_{nj}(z) f(z)| \to 0
\]
as $n\to \infty$, since $W_{n j}$ is a homogeneous polynomial of degree $n$.
It is also clear that $|f_{n j}(z)|\le \|f\|_{C(\bdd)}$ for all $z\in\bdd$. Therefore, by the dominated convergence theorem, we have
		$$
		\lim_{n\to\infty} \int_{\bd} |f_{n j}(z)|^q\, d\mu(z) = 0.
		$$

Using \eqref{e_left} and \eqref{e_right},
		we let $n$ tend to infinity in the inequality \eqref{RCM_sum} and we obtain
		\begin{multline*}
			\|f\|_{H^p} 
\le C\limsup_{n\to\infty} 
\sum_{j=1}^J \left( \int_{\bd} |f_{n j}(z)|^q\, d\mu(z) + \int_{\spd} |f_{n j}(z)|^q\, d\mu(z) \right)^{1/q}\\
\le C\limsup_{n\to\infty} \left(\int_{\spd} |f(z)|^q \, d\mu(z)\right)^{1/q} 
= C\|f\|_{L^q(\nu)}, 
		\end{multline*}
		as required.
	\end{proof}
	
	\subsection{Key lemma}
To characterize the $q$-reverse Carleson measures for $H^p$, $p\neq q$, 
we need Lemma~\ref{complex_analysis_goes_away} below. 
By this lemma, we can substitute arbitrary continuous functions into the inequality \eqref{q-RCM-Hardy-balayaged},
replacing the $H^p$-norm by $L^p$-norm.

	\begin{lem}\label{complex_analysis_goes_away}
		Suppose that $\mu$ is a $q$-reverse Carleson measure for $H^p$, where $p>0$ and $q>0$. 
Let $\nu=\mu|_{\spd}$, $\nu\in M_+(\spd)$. Then there exists a constant $C > 0$ such that
		\begin{equation}\label{plugging_cont_f}
			\|f\|_{L^p(\spd)}\leq C \|f\|_{L^q(\nu)}
		\end{equation}
		for every $f\in C(\spd)$.
	\end{lem}
	
	\begin{proof}
		It is enough to prove the lemma for the polynomials.
		Indeed, for any continuous function $f$, the inequality \eqref{plugging_cont_f} 
then follows by the uniform approximation.

Fix a polynomial $\ph$. The corresponding Hankel operator is defined by
\[
V_\ph [g] = \ph g - \KK [\ph g],\quad g\in H^2(\spd),
\]
where $\KK: L^2(\si) \to H^2(\spd)$ is the orthogonal projection.
Since $\ph$ is a smooth function, the operator $V_\ph: H^\infty (\spd) \to C(\spd)$ is compact 
(see, e.g., \cite[Theorem~6.5.4]{Ru80}). Therefore,
\begin{equation}\label{e131}
[
\|g_n\|_{C(\spd)} \le 1\ \textrm{and}\ g_n \to 0 \
\textrm{weakly in}\ H^2(\spd)
] \Rightarrow
\|V_\ph g_n\|_{C(\spd)} \to 0;
\end{equation}
see, for example, \cite[Equation~(5.1)]{DV26}.
Let $W_{n j}$ be those provided by Proposition~\ref{p_RW}.
We apply property~\eqref{e131} to $W_{n j}$ in the place of $g_n$.
Thus,
\begin{equation}\label{e_Hankel}
\|\ph W_{n j} - \KK [\ph W_{n j}]\|_{C(\spd)} \to 0
\end{equation}
as $n\to\infty$.

By Lemma~\ref{general_balayage}, we have
\[
\|\KK[\ph W_{n j}]\|_{H^p} \le \|\KK[\ph W_{n j}]\|_{L^q(\nu)}.
\]
Thus,
\begin{equation}\label{e_main}
\sum_{j=1}^J \|\KK[\ph W_{n j}]\|_{H^p} \le \sum_{j=1}^J \|\KK[\ph W_{n j}]\|_{L^q(\nu)}.
\end{equation}
Now, we estimate both sides of the above inequality.
On the one hand, if $n$ is sufficiently large, then, by \eqref{e_Hankel} and \eqref{e_RWlow},
\[
\begin{split}
\er\|\ph\|_{L^p(\si)} +
\sum_{j=1}^J \|\KK[\ph W_{n j}]\|_{H^p} 
&\ge C \sum_{j=1}^J \|\ph W_{n j}\|_{L^p} \\
&= C\sum_{j=1}^J \left( \int_{\spd} |\ph W_{n j}|^p\, d\si \right)^{\frac{1}{p}} \\
&\ge C \left( \int_{\spd} \sum_{j=1}^J |\ph|^p |W_{n j}|^p\, d\si \right)^{\frac{1}{p}} \\
&\ge C\delta \|\ph\|_{L^p(\si)}.
\end{split}
\]
Selecting $\er = C\delta/2$, we conclude that
\begin{equation}\label{e_low}
\sum_{j=1}^J \|\KK[\ph W_{n j}]\|_{H^p} \ge C \|\ph\|_{L^p(\si)}.
\end{equation}

On the other hand, for sufficiently large $n$, we have by \eqref{e_Hankel} and \eqref{e_RWup},
\begin{equation}\label{e_hi}
\sum_{j=1}^J \|\KK[\ph W_{n j}]\|_{L^q(\nu)}
\le C \sum_{j=1}^J \|\ph W_{n j}\|_{L^q(\nu)} + \er\|\ph\|_{L^q(\nu)}
\le C \|\ph\|_{L^q(\nu)}.
\end{equation}

Combining \eqref{e_main}, \eqref{e_low} and \eqref{e_hi}, we obtain
\[
\|\ph\|_{L^p(\si)} \le C \|\ph\|_{L^q(\nu)},
\]
as required.
	\end{proof}
	
	\begin{rem}\label{reformulation}
		If $\mu$ satisfies the hypothesis of Lemma~\ref{complex_analysis_goes_away}, 
		then clearly $\mu$ is a $q$-reverse Carleson measure for the space $H^p$. 
Hence, the initial problem is reduced to characterization of the measures $\nu\in M_+(\spd)$ which satisfy the inequality \eqref{plugging_cont_f} for every $f \in C(\spd)$.
	\end{rem}

	\subsection{The case $q < p$}
	
	\begin{theorem}
		If $0 < q < p$, then there are no $q$-reverse Carleson measures for the space $H^p$.
	\end{theorem}
	
	\begin{proof}
		Let $Q\subset \spd$ be defined by  \eqref{e_qball}. 
		We apply inequality \eqref{plugging_cont_f} to the characteristic function $f= \mathbf{1}_Q$. 
		While this function is not continuous, it is pointwise approximable by a sequence of uniformly bounded continuous functions, hence, the dominated convergence theorem applies.
		
		Inequality \eqref{plugging_cont_f} with $f= \mathbf{1}_Q$ yields
		$$
		\si(Q)^{1/p} \leq C \nu(Q)^{1/q},
		$$
		and hence $\nu(Q)\ge c \si(Q)^{q/p}$. 

Consider a maximal family of pairwise disjoint sets $Q_j \subset \spd$, $\si(Q_j) = \si(Q)$. 
Then we have at least $C/\si(Q)$ sets.
Hence, summing the corresponding inequalities for these sets, we obtain
\[
\si(Q)^{\frac{q}{p}- 1} \le C \nu(\spd).
\]
The above inequality holds for arbitrarily small $\si(Q) > 0$ and $\frac{q}{p} - 1 < 0$, 
thus, we deduce that $\nu(\spd)$ is infinite. This contradiction finishes the proof.
	\end{proof}
	
	\subsection{The case $q > p$}
	\begin{theorem}
		Let $\mu\in M_+(\bbd)$ and let $\nu$ denote the restriction of $\mu$ to the boundary of the unit ball: $\nu=\mu|_{\spd}$. Let $\nu^a$ be the absolutely continuous part of $\nu$. Denote by $\beta$ the Radon--Nikodym derivative of $\nu^a$ with respect to Lebesgue measure $\si$ on the sphere. 
		If $0 < p < q$, then $\mu$ is a $q$-reverse Carleson measure for the space $H^p(\bd)$ if and only if 
		\begin{equation}\label{cond_beta}
			\frac{1}{\beta}\in L^{\frac{p}{q-p}}(\si).
		\end{equation}
	\end{theorem}
	\begin{proof}
		The sufficiency of \eqref{cond_beta} follows from the H\"{o}lder inequality; see \cite{DTV26} for details.

The proof of the reverse implication is also similar to the argument in \cite{DTV26}.
For the reader's convenience, we give below certain details.

Suppose that $\mu$ is a $q$-reverse Carleson measure for the space $H^p$.  
		Applying Lemma~\ref{complex_analysis_goes_away}, we see that the property
		$$
		\|f\|_{L^p(\spd)}\leq C \|f\|_{L^q(\nu)}
		$$
		holds for all continuous functions $f$. Let $g\in C(\spd)$.
Using the above inequality with $f = |g|^{1/p}$, we obtain

		$$
		\Big| \int_{\spd} g\, d\si \Big|\le \int_{\spd} |g|\, d\si \le C \Big(\int_{\spd} |g|^{q/p}\, d\nu \Big)^{p/q}.
		$$
		Let $r = q/p > 1$. Then, for the functional
		$$
		\phi(g) =  \int_{\spd} g\, d\si,
		$$
		we have $|\phi(g)|\le C \|g\|_{L^r(\nu)}$ for all $g\in C(\spd)$.
Hence, $\phi$ extends to a continuous linear functional on  $L^r(\nu)$. 
Thus, the measure $\si$ is absolutely continuous with respect to $\nu$ and 
we have $\si=\alpha\, \nu$ with $\alpha\in L^{r'}(\nu)$. 
Therefore, we can write $\nu^a = \frac{1}{\alpha}\, \si$. Thus $\beta = 1/\alpha$. 
		Then the condition $\int_{\spd} \alpha^{r^\prime}\, d\nu < \infty$ rewrites as follows:
		$$
		\infty> \int_{\spd} \alpha^{r^\prime}\, d\nu^a = 
\int_{\spd}\alpha^{r^\prime-1}\, d\si = \int_{\spd} \beta^{1-r^\prime}\, d\si.
		$$
		Since $r^\prime = \frac{q}{q-p}$, we obtain the required condition \eqref{cond_beta}.
	\end{proof}

\bibliographystyle{plain}
\bibliography{Carleson3}

\end{document}